\documentclass [twoside]{article}

\newtheorem{thm}{Theorem}[section]
\newtheorem{prop}[thm]{Proposition}
\newtheorem{remark}{Remark} 
\newtheorem{dfn}{Definition} 

\def\bI{{\bf J}}

\def\b1{{\bf 1}}

\def\bT{{\bf T}}

\def\bP{{\bf P}}
\def\bQ{{\bf Q}}
\def\bR{{\bf R}}

\def\bN{{\bf N}}
\def\bZ{{\bf Z}}
\def\bC{{\bf C}}

\def\Blat{\mbox{\it \raise2pt\hbox{"}\kern-2pt H}}

\begin{document}
\begin{center}
{\center{\Large{\bf
Combinatorial aspects of the mixed Hodge structure
} }}

 \vspace{1pc}
{ \center{\large{ Susumu TANAB\'E }}} \noindent
\begin{center}
 \begin{minipage}[t]{10.2cm}
{\sc Abstract.} {\em This is a review article on the combinatorial
aspects of the  mixed Hodge structure of  a Milnor fibre of the
isolated hypersurface singularity. We give a purely combinatorial
method to compute spectral pairs of the singularity under the
assumption of  simplicial Newton boundary and non-degeneracy of
the germ.}
 \end{minipage} \hfill
\end{center}
 \vspace{1pc}
\end{center}

{ \center{\section{Introduction}} }

 The aim of this article is to give a survey on the
 combinatorial aspects of the MHS of the
 cohomology of the Milnor fibre defined by a single function germ with isolated
 singularity (hypersurface singularity) under the assumption of
 simplicial Newton boundary and non-degeneracy of the germ.

 In the case of a convenient germ $f,$ A.G.Kouchnirenko
 \cite{Kouch} established a formula of Milnor number $\mu(f) = dim H^{n-1}(X_t)$
 for the Milnor fibre $ X_t = \{ x \in \bC^n; |x| \leq\epsilon, f(x)=t\}$
for small enough $\epsilon$ and generic $t \not =0.$ Based on a
fundamental theory by J.H.M.Steenbrink \cite{Steen1},
 V.I. Danilov \cite{Dan} (almost simultaneously Anatoly N.Kirillov
 \cite{Kir1} also) has calculated the MHS $H^{p,q} (H^{n-1}(X_t))$
under the assumption that $f$ is non-degenerate and simplicial
(see Definition ~\ref{dfn32}).

Despite these remarkable results, their description of
$H^{n-1}(X_t)$ is not refined enough to study more advanced
question on the topology and the analysis on the Milnor fibre
$X_t.$ For example to calculate the Gauss-Manin system of the
fibre integrals $\int_{\gamma_j(t)} \omega_i,$ $\gamma_j(t) \in
H_{n-1}(X_t),$ $\omega_i \in H^{n-1}(X_t)$ we must know the
precise disposition of representatives $\omega_i \in H^{n-1}(X_t)$
with respect to the Newton diagram $\Gamma(f).$ Or, at least, to
describe the basis $\{\omega_1, \cdots, \omega_\mu\}$ in terms of
integer points on $\bR_+^n$ by means of combinatorics associated
to $\Gamma(f).$ This task has been carried by A.Douai
\cite{Douai1} for the case $n=2$ and non-degenerate $f$ to obtain
a concrete expression of the Gauss-Manin system on $ H^{1}(X_t).$
So far as it is known to me, the question of combinatorial
description of the $H^{p,q}(H^{n-1}(X_t))$ is still open. Here we
try to give an answer to this question (\S 1, {\bf Algorithm} ).

Quite recently, an algorithm to compute $H^{p,q}(H^{n-1}(X_t))$
together with the monodromy action on it has appeared (see
\cite{Schul01}). It is implemented in the computer algebra system
{\sf SINGULAR} in the library {\sf gaussman.lib}. Everybody who
wants to verify combinatorial statements on
$H^{p,q}(H^{n-1}(X_t))$ can achieve it in computing non-trivial
examples by means of this extremely useful tool.

{ \center{\section{Computational algorithm for the MHS of the
vanishing cohomology}} }

We describe here the mixed Hodge structure of the local
(vanishing) cohomology of the Milnor fibre. From combinatorial
point of view, the local structure is considered as a combination
of combinatorics treated in the global case i.e. MHS of the
cohomology of an affine algebraic variety in a torus teated by
\cite{DX1}.

Let us consider a germ $f(x) \in \bC[[x_1, \cdots,x_n]]$ that
defines an isolated singularity at $x=0.$ That is to say dimension
$ \mu(f)$ (Milnor number) of the Milnor ring $A(f)$ defined below
is finite:
$$  A(f) :=
\frac{\bC[[x_1, \cdots,x_n]]}{\bigl< \frac{\partial f}{\partial
x_1}, \cdots,  \frac{\partial f}{\partial x_n}\bigr>\bC[[x_1,
\cdots,x_n]]}. \leqno(1.1) $$ For a convex set
$$\Gamma_+(f)  := {\rm convex\;\; hull\;\; of} \{ \vec \alpha  +\bR_+^n
; \vec \alpha \in supp(f)\setminus \{0\} \}, \leqno(1.2) $$ we
define Newton boundary of the germ $f(x),$ $\Gamma(f) := $ union
of all closed compact faces of $\Gamma_+(f).$

We call a germ $f(x)$ convenient if it allows a decomposition as
follows,
$$f(x)= g(x)+R(x),$$
with $g(x)=\sum_{i=1}^n a_ix_i^{n_i},$ $\prod_{i=1}^n a_i\not =0,$
$n_i \geq 2$ for all $i\in [1,n]$ and $supp(R) \subset
\Gamma_+(g). $
\begin{dfn}
A germ $f(x)$ is called non-degenerate with respect to its Newton boundary
$\Gamma(f)$ if for every closed face $\tau\in \Gamma(f)$
the system of equations
$$ f^\tau(x)= x_1 \frac{\partial f^\tau}{\partial x_1}=
\cdots=  x_n \frac{\partial f^\tau}{\partial x_n}=0,$$ has no
common solutions in $\bT^n = (\bC^\times)^n.$
\label{dfn31}
\end{dfn}
This notion is similar to that of $\Delta(f)-$regular polynomial
defined in the global case, but it treats only $\tau\in
\Gamma(f).$ Let us denote by $\hat \tau$ the convex hull of $\tau
\cup \{ 0 \}.$ We define a ring $S_{\hat \tau}\subset$ $\bC[x_1,
\cdots, x_n]$ of the polynomial ring as follows:
$$S_{\hat \tau}:= \bC \oplus \bigoplus_{\frac{\vec \alpha}{k} \in \hat \tau,
\exists k \geq 1} \bC \cdot x^{\vec \alpha}. $$

Then the non-degeneracy of $f(x)$ is known to be equivalent to the
finite dimensionality of the ring
$$ A_{\tau} := \frac{S_{\hat \tau}}{\bigl<
x_1\frac{\partial f^\tau}{\partial x_1}, \cdots, x_n\frac{\partial
f^\tau}{\partial x_n}\bigr> S_{\hat \tau}}. \leqno(1.3)
$$
 Let us denote by $\Gamma_{-}(f)$ union of all segments
connecting $\alpha \in \Gamma(f)$ and $ \{ 0 \}$ or equivalently
$\Gamma_{-}(f) = \bigcup_{\tau \subset  \Gamma(f)}\hat \tau.$ Let
us denote by $V_k$ $k-$dimensional volume of disjoint sets (there
are $_nC_k$ such sets in total) $\Gamma_{-}(f)\cap $ $\{$
$k-$dimensional coordinate planes with $(n-k)$ zero coordinates
$\}.$

In this situation,
we have the following theorem on the Milnor number $\mu(f).$

\begin{thm}(\cite{Kouch})
Let $f(x)$ be a germ
convenient and non-degenerate with respect to $ \Gamma(f),$
then we have
$$
\mu(f) = n! V_n - (n-1)!  V_{n-1} + \cdots +(-1)^n.\leqno(1.4)
$$
\label{thm31}
\end{thm}

\begin{dfn}
We introduce the notion of simplicial Newton boundary which means
that  for each   $\tau \subset \Gamma(f)$ the following inequality
holds
$$ \sharp \{ \Gamma_i \; {\rm face \;of \;} \Gamma(f); dim \;\Gamma_i = dim \; \tau +1, \tau \subset
\Gamma_i\} \leq n- dim \; \tau. $$
\label{dfn32}
\end{dfn}
As a matter of fact, we can formulate the above theorem by
Kouchnirenko in a more precise form. We introduce a new  $\bC-$
vector space $V_\tau$ associated to a face $\tau \in \Gamma(f)$
not contained in a coordinate plane.
$$V_\tau = A_\tau \setminus(\oplus_{\tau^{(1)} \in \tau}  A_{\tau^{(1)}}
\setminus(\oplus_{\tau^{(2)} \in \tau}  A_{\tau^{(2)}} \setminus(
\cdots \setminus \{0\})\cdots),$$ where $\tau^{(j)} \in \tau$
denotes a codimension $j$ face of $\tau$ contained in a coordinate
plane. Here we remark that though $\tau$ not contained in a
coordinate plane $\tau^{(j)} \in \tau,$ $j\in [1, dim\tau]$ may be
contained in a coordinate plane. We introduce another $\bC -$
vector space $W_\tau$ corresponding to the interior points of
$supp(V_\tau)$,
$$W_\tau = A_\tau \setminus(\oplus_{\tau^{(1)} \in \tau}  A_{\tau^{(1)}}
\setminus(\oplus_{\tau^{(2)} \in \tau}  A_{\tau^{(2)}} \setminus(
\cdots \setminus \{0\})\cdots),$$ where $\tau^{(j)} \in \tau$
denotes a codimension $j$ face of $\tau$ not necessarily contained
in a coordinate plane.

We say that a set $c(\sigma)$ is a copy of set $\sigma$ if the
relation $c(\sigma)= \pm \sigma + \vec w,$ for some $\vec w \in
\bZ^n$ holds. Further on we use the notation $c^{j}(\sigma),$
$j=1,2,\cdots$ to distinguish different copies of a set $\sigma.$

Proposition 2.6 of \cite{Kouch}, (5.6), (5.7) of \cite{Steen1}
entail the following.

\begin{prop}
1)  For $A_\tau,$ we have the following relations,
$$dim \; A_\tau=  (dim \;\tau +1)!
vol (\hat \tau).$$

2)  $$\mu(f) = \sum_{\tau \subset {\rm coordinate \; planes}}
(-1)^{n-1-dim \tau} dim\; A_\tau $$

3) $$A(f) \simeq \oplus_{(n-1){\rm dimensional \;faces\;}\tau
\subset \Gamma(f)}V_\tau. \leqno (1.5)
$$
 In case of repetitive appearances of $A_\gamma$'s, for some
face $\gamma$ in different $V_{\tau_1},$ $\cdots,$ $V_{\tau_k}$
($\gamma \subset \tau_1 \cap \cdots \cap \tau_k$) these  copies of
$A_\gamma$ (or rather $supp(A_\gamma)$) shall be shifted and
located anew in a way that they form a symmetry with respect to
the Hodge filtration of $A_{\tau_i}$ for some $i \subset [1,k].$

4) Let us denote by $s^{(\ell)}(\sigma) $ a shift of a set $\sigma
\in F^i/F^{i+1}$ to another properly chosen copy
$s^{(\ell)}(\sigma) \in F^{i-\ell}/F^{i-\ell+1}.$ Then we have
another representation as follows,
$$A(f) \simeq
\bigoplus_{\sigma \subset \Gamma(f)} \bigoplus_{\tau \subset
\sigma} \bigoplus^{n- dim \sigma -1}_{\ell=0} \bigoplus^{\;_{n-
dim \sigma -1} C_{\ell }}_{j=1} \bigoplus_{i=0}^{dim \tau}(-1)^{n-
dim \sigma -1 -\ell}
c^{j}\big(s^{(\ell)}\bigl(F^i/F^{i+1}(W_\tau)\bigr) \big).
\leqno(1.5)'
$$ Here different copies of $c^{j}\big(s^{(\ell)}(W_\tau) \big)$
shall be distributed in $\bigoplus_{\sigma \subset \Gamma(f)}
\bigoplus_{\tau \subset \sigma} F^{i-\ell}/F^{i-\ell+1}(A_{\tau}
),$ in such a way that $c^{j}\big(s^{(\ell)}(W_\tau) \big)$ $\cap$
$c^{j'}\big(s^{(\ell)}(W_\tau) \big) = \emptyset $ for all pairs
$j\not = j'. $
 \label{prop62}
\end{prop}
A precise way to arrange copies in accordance with the Hodge
filtration shall be explained in the {\bf Algorithm} below.

Further we shall establish a connexion between the volume of a
polyhedron and a set of integer points. Let $\tau $ be a
$(k-1)-$dimensional face of $\Gamma(f)$ and $\hat \tau$ be a
$k-$dimensional convex polyhedron. Let us denote by $\vec m_1,
\cdots, \vec m_r$ vertices of $\hat \tau \setminus\{0\}.$ We
consider the cone
$$
cone(\tau)= \{ \sum_{i=1}^r b_i \vec m_i;b_i \geq 0 \},
\leqno(1.6)
$$ associated to $\tau.$
We introduce a grading on the algebra $S_{\hat \tau}.$ First we
consider a piecewise linear function $h:\bN^n \rightarrow \bN$
satisfying $h|_{\Gamma(f)}=1.$ Then there exists $M>0$ such that
$h(\alpha) \subset \frac{1}{M}\bN$ for all $\alpha \in \bN^n.$ We
define $\phi= M\cdot h|_{\bN^n}.$ Let us denote by ${\mathcal
A}_q$ algebra of polynomials written as a linear combination of
monomials $x^\alpha, \phi(\alpha) \geq 1.$ Denote by ${\mathcal
A}_q(\tau)$ subalgebra of  polynomials of ${\mathcal A}_q$ whose
supports are contained in $cone(\tau).$ Then we can consider the
Poincar\'e polynomial of $S_{\hat \tau}$ defined by
$$P_{S_{\hat \tau}}(t):=
\sum_{q=0}^\infty dim_{\bC} \big({\mathcal A}_q(\tau)/ {\mathcal
A}_{q+1}(\tau) \big)t^q.$$ Then we have the following relationship
$$
k! vol_k(\hat \tau) = \sharp \{\bZ^n \cap \{cone(\tau) \setminus
\bigcup_{i=1}^r (\vec m_i +cone(\tau)) \} \} = P_{S_{\hat
\tau}}(t)(1-t)^k|_{t=1}. \leqno(1.7)
$$
Here we recall the fundamental theorem from \cite{Steen1} (3.10).
To formulate it, we need to introduce preparatory notions. Let us
consider a resolution of singularity $X_0,$ that is to say a
proper mapping $\rho: Y \rightarrow \bC^n$ from a smooth algebraic
variety $Y \supset \bC^n$ such that 1) $\rho$ is an isomorphism on
$\bC^n \setminus\{0\}$ and 2)$E=\rho^{-1}(X_0)$ is a divisor on
$Y$ with transversal intersections. Let $E_0$ be the proper image
of $X_0$ through $\rho,$ i.e. the closure of $\rho^{-1}(X_0
\setminus \{0\})$ in $Y.$ Let us denote by $E_1, \cdots E_N$ the
remaining irreducible components of $E.$ Assume that $E= E_0 +
\sum_{i=1}^N m_iE_i$ with multiplicities $m_i$ of the divisor
$E_i.$ Let $M$ be the least common multiplier (l.c.m.) of $m_1,
\cdots, m_N.$ We consider a covering $\pi:\tilde \bC \rightarrow
\bC$ that sends $z$ to $z^M.$ For the pair of mappings $(f,\pi)$
we denote the fibre product $Y \times_{\bC} \tilde \bC$ by $\tilde
X. $ Let $D_i = \pi^{-1}(E_i)_{red}, i \in [1,N]$ be the reduced
part of $\pi^{-1}(E_i).$ If we consider the morphism $\tilde f :
\tilde X \rightarrow \tilde \bC,$ and its special fibre $D:=
\tilde f^{-1}(0), $ then we have $D =\sum_{i=1}^m D_i.$ We will
use the notations,
$$D^{(k)}= \coprod_{i_0 < \cdots <i_k} \big(D_{i_0}\cap \cdots \cap D_{i_k}
\big)_{red}, 'D^{(r)}= \coprod_{0<i_0 < \cdots <i_r}
\big(D_{i_0}\cap \cdots \cap D_{i_r} \big)_{red}.$$ Under these
circumstances we have the following theorem (\cite{Dan},
\cite{Steen1}) on the vanishing cohomology $H^{r+k}(X_\infty)$.
\begin{thm}
There exists a spectral sequence $E_1^{r,k}$ converging to
$H^{r+k}(X_\infty)$ satisfying the following properties.

1) It converges to the weight filtration on $H^{r+k}(X_\infty),$
i.e. $E_\infty^{r,k}= Gr^W_k \;H^{r+k}(X_\infty),$

2) It degenerates at the term $E_2$ and $E_2 = E_\infty.$

3) The $E_1$ term is given by the formulae
$$
\begin {array}{ccccccc}
E_1^{r,k}=& \oplus_{i \geq 0} H^{k+2r-2i}(D^{(2i-r)})(r-i)& for&
r<0, &
\\
=&H^{k}('D^{(r)}) \bigoplus \big(\oplus_{i >r }
H^{k+2r-2i}(D^{(2i-r)})(r-i) \big)& for& r\geq 0 & & .&
\\
\end{array}$$
\label{thm33}
\end{thm}

We can classify the elements of $A_{\tau}$ after their eigenvalues
under the action $x \rightarrow \zeta_\ast(x)=
\zeta^{-h(\alpha)}x^\alpha$ with $\zeta=e^{\frac {2 \pi \sqrt
-1}{M}}$ that coincides with the action $T_s$ of the semisimple
part of the monodromy $T=T_s \cdot T_u,$ where $T_u$ denotes the
unipotent part of $T.$

Let us introduce the Poincar\'e polynomial of ${\mathcal
A}_q(\tau)/ {\mathcal A}_{q+1}(\tau)$ in taking  the monodromy
action $\zeta_\ast$ into account,
$$ P_{{\mathcal A}_q(\tau)/
{\mathcal A}_{q+1}(\tau)}(t) := \sum_{0 < \chi < 1} h_\chi^{q, dim
\tau-q} t^\chi. \leqno(1.8)$$
$$ \tilde P_{{\mathcal A}_q(\tau)/
{\mathcal A}_{q+1}(\tau)}(t) := h_1^{q,q}t^q. \leqno(1.9)$$ where
$$h_\chi^{q, dim \tau-q}:= \sharp \{x^\alpha \in {\mathcal A}_q(\tau)/ {\mathcal
A}_{q+1}(\tau); h(\alpha)= \chi + q \},$$
$$h_1^{q,q}:=\sharp \{x^\alpha \in {\mathcal
A}_q(\tau)/ {\mathcal A}_{q+1}(\tau); h(\alpha)= q\}.$$

The main theorem of $\cite{Dan}$ can be formulated as follows,
\begin{thm}We suppose that $\Gamma(f)$ is a simplicial Newton
boundary. Then Poincar\'e polynomials $(1.8),$ $(1.9)$ satisfy the
following relations,
$$P_{{\mathcal A}_q(\tau)/
{\mathcal A}_{q+1}(\tau)}(t) =(-1)^{dim \tau -q}\sum_{{\rm all
\;faces}\gamma \subset \tau} \sum_{k \geq 0}(-1)^k _{dim\; \gamma
+1}C_{p+k+1} \big( \sum_{\alpha \in (k+1)\hat \gamma}t^{h(\alpha)}
- \sum_{\alpha \in k\hat \gamma}t^{h(\alpha)} -  \sum_{\alpha \in
k \gamma}t^{h(\alpha)}\big), \leqno(1.10)$$
$$ \sum_{q \geq 0}\tilde P_{{\mathcal A}_q(\tau)/
{\mathcal A}_{q+1}(\tau)}(t) = \sum_{{\rm all \;faces}\gamma
\subset \tau} (t-1)^{dim \;\gamma}. \leqno(1.11)$$ \label{thm34}
\end{thm}

Let us recall fundamental notions around the spectral pairs of the
singularity  that reflect the interplay between the monodromy
action $T$ and the MHS of $H^{n-1}(X_\infty)$ \cite{Steen2}. The
MHS on $H^{n-1}(X_\infty)$ consists of an increasing weight
filtration $W_\cdot$ and a decreasing Hodge filtration $F^\cdot$ (
\cite{Steen1}). Let $T_s$ be the semisimple part of $T,$ and $T_u$
unipotent, then $T_s$ preserves the filtration $F^\cdot$ and
$W_\cdot$ whereas $N= log\; T_u$ satisfies $N(W_i) \subset
W_{i-2}$ and $N(F^p) \subset F^{p-1}.$ For eigenvalue $\chi$ of
$T,$ we define
$$ H^{p,q}_\chi := Ker \big(T_s -\chi \cdot id_\mu ; Gr^W_{p+q} \tilde H^{n-1}(X_\infty)
\big),$$
$$ dim\;H^{p,q}_\chi  = h^{p,q}_\chi,$$
where $\tilde H^{n-1}(X_\infty)$ denotes the reduced cohomology,
$Gr^W_i = W_i/W_{i-1},$ and $Gr^p_F = F^p/F^{p+1}.$ For $\alpha
\in \bQ$ and $w \in \bZ$ we define integers $m_{\alpha,w}$ as
follows. Write $\alpha = n-1-p-\beta$ with $0 \leq \beta <1$ and
let $\chi = e^{-2\pi \sqrt -1 \alpha}.$ If $\chi \not =1$ then
$m_{\alpha,w}=h^{p,w-p}_\chi$ while
$m_{\alpha,w}=h^{p,w+1-p}_{\chi=1}.$ The spectral pairs are
collected in the invariant
$$ Spp(f)= \sum m_{\alpha,w}(\alpha,w), \leqno(1.12)$$
to be considered as an element of the free abelian group on $\bQ
\times \bZ.$ It is known that $Spp(f)$ is invariant under the
symmetry $( \alpha,w ) \rightarrow (n-2-\alpha, 2n-2-w)$
\cite{Steen2}, Theorem 1.1, (ii).

 Theorem ~\ref{thm34} entails the
relations
$$\sum_{dim \hat \tau
=d}P_{{\mathcal A}_q(\tau)/ {\mathcal A}_{q+1}(\tau)}(t)= \sum_{0
<\chi <1} h_\chi^{q, d-q} t^\chi, \leqno(1.13)$$
$$ \sum_{q \geq 0}\sum_{\tau \subset \Gamma(f)} \tilde P_{{\mathcal A}_q(\tau)/
{\mathcal A}_{q+1}(\tau)}(t) = \sum_{q \geq 0}h^{q,q}_1t^q.
\leqno(1.14)$$ As a corollary we have,
$$ h_{\chi}^{n-1-p, n-1-q}= h_{\chi^{-1}}^{p, q},\; h_{\chi=1}^{p,p}= h_{\chi=1}^{n-p, n-p}.
\leqno(1.15)$$ We can write down the formula $(1.10)$ in a more
combinatorially clear way,
$$ h^{p, dim\; \tau -p}_{\chi\not =1}(D_{\hat \tau})
= (-1)^{dim\; \tau -p} \sum_{{\rm all \;faces}\gamma \subset \tau}
\sum_{k \geq 0}(-1)^k _{dim\; \gamma +1}C_{p+k+1} \big(
\ell^\ast((k+1)\hat \gamma) - \ell^\ast( k\hat \gamma) -
\ell^\ast( k \gamma)\big) \big), \leqno(1.16)$$ where $D_{\hat
\tau}= \bP_{\tilde \tau} \cap \tilde X$ for $\tilde \tau$
suspension of $(\tau,0) \subset \bR^{n+1}$ with $(0,\cdots,0,M)
\in \bR^{n+1}.$

{\bf Algorithm}

This is the unique original part of this article.

Further we give an algorithm to get a basis of $A(f)$   in a
purely combinatorial way under the assumption that $\Gamma(f)$ is
a simplicial boundary and $f$ is a non-degenerate germ. We shall
achieve this task in making the decomposition of $A(f)$ in $(1.5)$
more precise. Though $(1.5)'$ gives us a more detailed description
of   $A(f)$ than that of $(1.5)$ it turns out less convenient for
the construction of an algorithm.  We remark that $(1.5)'$ has
been obtained from the evident relation $V_\sigma = \bigcup_{\tau
\subset \sigma} W_{\tau} $ or $P_{A_\sigma}(t) =(1-t)^{dim\;
\sigma +1}P_{S_{\hat \sigma}}(t) = \sum_{\tau \subset \sigma}
P_{W_{\tau}}(t) $ on the level of Poincar\'e polynomials. Thus
$(1.5)'$ contains exactly the same combinatorial
 informations as in $(1.5).$

 Let $\vec m_1, \cdots, \vec m_k$ be vertices of a
$k-$dimensional simplex face $\tau$ (if necessary we divide a
non-simplex face into a sum of simplices). Here we remark the fact
that for two simplices $\tau_1, \tau_2$ whose sum give a face
$\Delta \subset \Gamma(f)$ i.e. $\Delta = \tau_1\cup \tau_2$ and
whose intersection is again a simplex $\gamma$ ;$\gamma = \tau_1
\cap \tau_2, $ we have
$$ P_{S_{\hat \Delta}}(t) = P_{S_{\hat \tau_1}}(t) + P_{S_{\hat \tau_1}}(t) -
P_{S_{\hat \gamma}}(t).$$ Thus the following procedure has
meaning.
\begin{dfn} Simplex subdivision $\delta_1, \cdots, \delta_m$
of faces of $\Gamma(f)$ means that for each $(n-1)$ dimensional
compact face $\gamma \subset \Gamma(f),$ there exists a
subdivision of it into a sum of $(n-1)-$ dimensional simplices,
$$\gamma = \bigcup_{i \in I(\gamma)}\delta_i,$$ for a set of
indices $I(\gamma) \subset [1, \cdots, m]$ associated to $\gamma.$
Consequently,
$$\Gamma_{-}(f) = \bigcup_{i =1}^m \hat \delta_i,$$
is a subdivision into $n$ dimensional simplices $\hat \delta_i,$
$1 \leq i \leq m.$
 \label{dfn}
\end{dfn}

We describe a combinatorial algorithm (not unique) to get a basis
of $A(f)$ consisting of several steps.

1) {\it For a $(n-1)$dimensional simplex $\tau$ (whose vertices
are $\vec v_1,$ $\cdots,$ $\vec v_n$) of a simplex subdivision, we
construct the parallelepiped
$$ B_\tau:= \{\bR^n \cap \{cone(\tau) \setminus
\bigcup_{i=1}^n (\vec v_i +cone(\tau)) \} \}. \leqno(1.17)$$ The
inclusion relation $B_\tau \supset supp (A_\tau)$ $\supset supp(
V_\tau)$ can be easily seen from $(1.6).$ For fixed subset of
indices $\bI \subset \{1, \cdots, n\}$ each vertex of the
parallelepiped has the form
$$ \vec v(\bI) := \sum_{i \in \bI} \vec v_i,$$
where no repetition of indices is allowed.}

2) {\it To consider the set $G_\tau =$ $B_\tau$ $\setminus$ $\{$
all open skeletons of dimension less than $(n-1)$ contained in
$F^0/F^1(A_{\tau})$ $\}.$ In other words $G_\tau = supp(W_\tau)$}.

As a special case of copy, we introduce the notion of canonical
copy $c_\tau(\alpha)$ of a point $\alpha$ with respect to a
$(n-1)$dimensional simplex $\tau$ of the simplex subdivision
(whose vertices are $\vec v_1,$ $\cdots,$ $\vec v_n$) that means
the points $\alpha,$ $c_{\tau}(\alpha)$ are symmetrically located
with respect to $\frac{1}{2} \sum_{i=1}^n \vec v_i,$
$$ c_{\tau }(\alpha)+\alpha =\sum_{i=1}^n \vec v_i.   \leqno(1.18)$$
We shall choose basis of $A(f)$ in such a way that the symmetry
property of Hodge numbers $(1.15)$ can be realized. As for the
integer points of $A_\tau$ on the intermediate Hodge filtration
level $F^i/F^{i+1}(A_{\tau}),$ $1 \leq i \leq n-2,$ the points of
$G_\tau$ already realize this symmetry property. This can be seen
from the arguments of \cite{DX1}, \S 5 where essentially
$supp(A_\tau)$ is combinatorially described. Moreover, for
different $(n-1)-$simplices $\tau_1$ and $\tau_2$ from simplex
subdivision, $ F^i/F^{i+1}(A_{\tau_1})$ $\cap$
$F^i/F^{i+1}(A_{\tau_2})= \emptyset, $ $1 \leq i \leq n-2.$ This
can be seen from the fact that  $F^i/F^{i+1}(A_{\tau_j}) \subset
cone(\tau_j), j=1,2.$ Thus we shall further first care about the
choice of $supp(A_\tau)$ on the extremal Hodge filtration levels
$F^0/F^{1}(A_{\tau})$ and $F^{n-1}/F^{n}(A_{\tau}).$

3) {\it To count the number of interior points of each canonical
copy $c_{\tau}(\hat \tau^{int})$ of $\hat \tau^{int}$ in $
G_\tau$, located on the Hodge filtration level
$F^0/F^{1}(A_{\tau}).$}

4) {\it For every $(n-1)$ simplex $\tau$ from simplex subdivision
to exclude faces from $G_\tau$, contained in
$F^{n-1}/F^n(A_{\tau}),$ that are located on some coordinate
plane.}


The following two measures 5), 6) are to be taken to cope with
repetitive appearances of  $A_\gamma$'s mentioned in the
Proposition ~\ref{prop62}, 3).

 5) {\it Suppose that $\Delta_1,$
$\cdots, \Delta_k$ are $(n-1)$ simplices from a simplex
subdivision of faces of $\Gamma(f)$ such that $\hat \Delta_1 \cap
\cdots \cap \hat \Delta_k \not = \emptyset.$ To choose a canonical
copy $c_{\Delta_i}(\sigma^{int})$ of each open skeleton
$\sigma^{int}$ of   $\hat \Delta_1 \cap \cdots \cap \hat \Delta_k$
with respect to a simplex $ \Delta_i$ that is to be chosen in
dependence of $\sigma^{int}.$ If the open skeleton $\sigma^{int}$
has another expression like $\sigma^{int} \subset \hat \gamma_1
\cap \cdots \cap \hat \gamma_{k'}$ for another pair of simplices
of a simplex subdivision $\{ \Delta_1, \cdots, \Delta_k\}$ $\not
=$ $\{ \gamma_1,\cdots, \gamma_{k'} \},$ we do not add any of
canonical copies $c_{\gamma_j}(\sigma^{int}),\; j \in [1, k'].$}

This procedure is necessary to recover these integer points that
are located on the intersection $\hat \Delta_1 \cap \cdots \cap
\hat \Delta_k$ on the level of $F^0/F^1(A_{\Delta_i})$ for some
unique $i \in [1,k].$

For example, in $f_3$ case below (see {\bf 2.3}) $(1,1, 1)\in (0,
\vec v_0 )^{int}$ contained in $\hat \Gamma_1 \cap  \hat
\Gamma_{2}, $ $\hat \Gamma_2 \cap  \hat \Gamma_{3}$ and $\hat
\Gamma_3 \cap  \hat \Gamma_{1}.$ If we add the canonical copy
$c_{\Gamma_2}((0, \vec v_0 )^{int}) = (\vec v_1 + \vec v_2, \vec
v_0 + \vec v_1 + \vec v_2 )^{int}$   to $G_{\Gamma_2}$  neither
$c_{\Gamma_1}((0, \vec v_0 )^{int})$ nor $c_{\Gamma_3}((0, \vec
v_0 )^{int})$ is needed any more.

6) {\it Furthermore if $dim(\hat \Delta_1 \cap \cdots \cap \hat
\Delta_k)=$ $dim \;\sigma^{int}$ we shall add other not canonical
copies $c^{2}(\sigma^{int}),$ $\cdots,$ $c^{k-1}(\sigma^{int})$
(in understanding $c^{1}(\sigma^{int})=\sigma^{int},$
$c^{k}(\sigma^{int})=c_{\Delta_i}(\sigma^{int})$ of the procedure
5) above) such that
$$ c^{[\frac{k+1}{2}]+j}(\sigma^{int})\in
F^{[\frac{n+1}{2}]+j}/F^{[\frac{n+1}{2}]+j+1}(A_{\Delta_{i_j}})\leqno(1.19)$$
for $2-[\frac{k+1}{2}]\leq j \leq k-[\frac{k+1}{2}]-1.$ Thus they
produce a symmetry with respect to the Hodge filtration
$F^\cdot.$}

In the case of simplicial Newton boundary $\Gamma(f)$ we have $k
\leq n$ thus the above procedure can be realized so that $(1.19)$
holds in such a way that $c^{j}(\sigma^{int}) \in
G_{\Delta_{i_j}}$ and $ \Delta_{i_j} \not = \Delta_{i_{j'}} $ for
all pairs $j \not = j'.$ On the contrary, if $\Gamma(f)$ is not
simplicial, such a simple construction is already impossible. This
situation explains why Danilov restricted himself to the
simplicial Newton boundary case in \cite{Dan}.

For example, see $(2.1.1)$, $(2.1.2)$ and $(2.1.3)$ below.

7) {\it Add zero dimensional faces (i.e. vertices) of $\Delta_j$
not belonging to the coordinate plane and their canonical copies
with respect to $\Delta_j$ only once for each.}

Making use of the above basis, one can calculate the MHS of
$A(f).$

8) {\it We classify all points from $x^{\vec \alpha} \in A(f)$
according to their position with respect to faces of simplex
subdivision $\delta_1, \cdots, \delta_m.$  That is to say to find
$\delta_i$ such that $$ x^{\vec \alpha+\b1} \in B_{\bar
\delta_i},$$ where $x^\b1 = x_1\cdots x_n.$
 }

9) {\it To evaluate $ h(\vec \alpha+\b1) $ by means of the
piecewise linear function $h$ such that $h|_{\delta_i}=1$
introduced just after $(1.6).$}

Further on, in contrast to the above procedures where we meant by
$\sigma$ an intersection of $(n-1)-$dimensional simplices of a
simplex subdivision, we consider as $\sigma$ only faces of
$\Gamma(f)$ whose dimension vary from $0$ to $n-1.$

 10) {\it (${\chi \not= 1}$ case )   If
$h(\vec\alpha+\b1)=n-1-\beta -p $ for $0 \leq p \leq n-1,$  $0
<\beta<1$, then $ x^{\vec \alpha} \in H_{\chi \not =1}^{p,q}.$
Here the index $q$ can be chosen in the following way. For $p <[
\frac{n-1}{2}]$ the index $q$ is to be chosen $q=dim\;
\sigma^{int} -1 >0$ if $\vec \alpha+\b1$ belongs to one of the
copies of $ \sigma^{int}.$ While for  $p >[ \frac{n-1}{2}]$ the
index $q$ is to be chosen $q=n- dim\; \sigma^{int} \geq 0$ under a
parallel situation. All other cases ($h^{p,q}_{\chi \not =1}$ )
except $H^{\frac{n-1}{2},\frac{n-1}{2}}_{\chi \not =1},$ ($n:odd$)
can be recovered from the above data making use of the relation
$(1.15)$ $h^{p,q}_{\chi}= h^{n-1-p,n-1-q}_{\chi^{-1}}$ realized by
taking proper copies. The exceptional case has a following
expression,
$$ H^{\frac{n-1}{2},\frac{n-1}{2}}_{\chi \not =1} \cong \{ \vec \alpha+\b1 \in \bigcup_{i=1}^m
B_{\delta_i}^{int}; \frac{n-1}{2} < h (\vec \alpha+\b1) <
\frac{n+1}{2} \}.$$ }

11) {\it (${\chi = 1}$ case ) If   $ h(\vec \alpha+\b1)=n-1 -p $
for $0 \leq p \leq n-1$ and $\vec \alpha+\b1$ belongs to one of
the copies of $ \sigma^{int}$, then $ x^{\vec \alpha} \in H_{\chi
=1}^{p+1,q}.$ Here the index $q$ can be chosen as $q=dim\;
\sigma^{int}>0$ if $\vec \alpha+\b1$ belongs to one of the copies
of $ \sigma^{int},$ while for  $p >[ \frac{n-1}{2}]$ the index $q$
is to be chosen $q=n- dim\; \sigma^{int}
>0$ under a parallel situation. All cases can be recovered from the above data
making use of the relation $(1.15)$ $h^{p,q}_{\chi=1}=
h^{n-p,n-q}_{\chi=1}$ realized by taking proper copies.}

\begin{remark}
The choice of the representative $mod {J_{f,\Delta}} $ in
$B_\Delta$ does effect  not only on the weight filtration but also
on the Hodge filtration (see examples below). \label{remark61}
\end{remark}

{ \center{\section{Examples}} }

We show examples of calculus by means of the computer algebra
system for computation $\sf SINGULAR.$ One can find an
introduction to algorithms to compute monodromy related invariants
(namely spectral pairs) of isolated hypersurface singularities in
\cite{Schul1}. In the sequence, we use the notation $[i]'=[i]xy$
for $\bf 2.1$ and $[i]'=[i]xyz$ for $\bf 2.2,$ $\bf 2.3.$   In the
description of the spectral pairs we use the convention
$((\alpha,w),m_{\alpha,w})$ under the notation of $(1.12).$ We see
that the rational monodromy $\alpha_i$ of the basis $[i]$ is
expressed as $\alpha_i = h([i]')-1$ for piecewise linear function
$h(\cdot)$ introduced just after $(1.6).$

 {\bf 2.1} Let us begin with a polynomial in two variables,
$$f_1=x^{15}+x^6y^4+x^3y^6+y^{12}.$$
Here and further on, we shall make use of the notational
convention $xiyjzk= x^iy^jz^k.$ The algebra $A(f_1)$ (rank
$A(f_1)= 94$) has the following basis, $[1]=xy13 , [2]=y13 ,
[3]=xy12 , [4]=y12 , [5]=xy11 , [6]=y11,$ $ [7]=xy10 , [8]=y10 ,
[9]=xy9 , [10]=y9 , [11]=xy8 , [12]=y8 ,$ $[13]=xy7 , [14]=y7 ,
[15]=xy6 , [16]=y6 , [17]=x2y5 , [18]=xy5 ,$ $[19]=y5 , [20]=x5y4
, [21]=x4y4 , [22]=x3y4 , [23]=x2y4 , [24]=xy4 ,$ $[25]=y4 ,
[26]=x8y3 ,$ $[27]=x7y3 , [28]=x6y3 , [29]=x5y3 ,$ $ [30]=x4y3 ,
[31]=x3y3 $, $[32]=x2y3 , [33]=xy3 , [34]=y3 ,$ $[35]=x16y2 ,
[36]=x15y2 ,$ $[37]=x14y2 , [38]=x13y2 , [39]=x12y2 ,$ $[40]$ =
$x11y2 , [41]$ = $x10y2 , [42]$ = $x9y2 , [43]$ = $x8y2 , [44]$ =
$x7y2 , [45]$ = $x6y2 , [46]$ = $x5y2 , [47]$ = $x4y2 , [48]$ =
$x3y2 , [49]$ = $x2y2 , [50]$ = $xy2 , [51]$ = $y2 , [52]$ = $x19y
, [53]$ = $x18y , [54]$ = $x17y , [55]$ = $x16y , [56]$ = $x15y ,
[57]$ = $x14y , [58]$ = $x13y , [59]$ = $x12y , [60]$ = $x11y ,
[61]$ = $x10y , [62]$ = $x9y , [63]$ = $x8y , [64]$ = $x7y , [65]$
= $x6y , [66]$ = $x5y , [67]$ = $x4y , [68]$ = $x3y , [69]$ = $x2y
, [70]$ = $xy ,$ [71]$ = $y , [72]$ = $x22 , [73]$ = $x21 , [74]$
= $x20 , [75]$ = $x19 , [76]$ = $x18 , [77]$ = $x17 , [78]$ = $x16
, [79]$ = $x15 , [80]$ = $x14 , [81]$ = $x13 , [82]$ = $x12 ,
[83]$ = $x11 , [84]$ = $x10 , [85]$ = $x9 , [86]$ = $x8 , [87]$ =
$x7 , [88]$ = $x6 , [89]$ = $x5 , [90]$ = $x4 , [91]$ = $x3 ,
[92]$ = $x2 ,$ [93]$ = $x , [94]$ = $1.$

The spectral pairs are calculated as follows,
\noindent

$ ($$ ($-19/24,1$ )$,1$ )$,$ ($$ ($-43/60,1$ )$,1$ )$,$ ($$
($-2/3,2$ )$,1$ )$,$ ($$ ($-13/20,1$ )$,1$ )$,$ ($ $ ($-7/12,1$
)$,3$ )$ , $ ($ $ ($-31/60  , 1$ )$ , 1$ )$ , $ ($ $ ($-1/2 , 2$
)$ , 1$ )$ , $ ($ $ ($-1/2 , 1$ )$ , 1$ )$ , $ ($ $ ($-11/24 , 1$
)$ , 1$ )$ , $ ($ $ ($-9/20 , 1$ )$ , 1$ )$ , $ ($ $ ($-13/30 , 1$
)$ , 1$ )$ , $ ($ $ ($- 5/12 , 1$ )$ , 1$ )$ , $ ($ $ ($-23/60 ,
1$ )$ , 1$ )$ , $ ($ $ ($-3/8 , 1$ )$ , 1$ )$ , $ ($ $ ($-11/30 ,
1$ )$ , 1$ )$ , $ ($ $ ($-1/3 , 2$ )$ , 1$ )$ , $ ($ $ ($-1/3 , 1$
)$ , 1$ )$ , $ ($ $ ($-19/60,1$ )$,1$ )$,$ ($ $ ($-3/10,1$ )$,1$
)$,$ ($ $ ($-7/24,1$ )$,1$ )$,$ ($ $ ($-1/4,1$ )$ , 4$ )$ , $ ($ $
($-7/30 , 1$ )$ , 1$ )$ , $ ($ $ ($-13/60 , 1 $ )$ , 1$ )$ , $ ($
$ ($-11/60 , 1$ )$ , 1$ )$ , $ ($ $ ($-1/6 , 1$ )$ , 4$ )$ , $ ($
$ ($-3/20 , 1$ )$ , 1$ )$ , $ ($ $ ($-1/8 , 1$ )$ , 1$ )$ , $ ($ $
($-7/60 , 1$ )$ , 1$ )$ , $ ($ $ ($-1/1 0 , 1$ )$ , 1$ )$ , $ ($ $
($-1/12 , 1$ )$ , 4$ )$ , $ ($ $ ($-1/20 , 1$ )$ , 1$ )$ , $ ($ $
($-1/24 , 1$ )$ , 1$ )$ , $ ($ $ ($-1/30 , 1$ )$ , 1$ )$ , $ ($ $
($-1/60 , 1$ )$ , 1$ )$ , $ ($ $ ($ 0 , 1$ )$ , 4$ )$ , $ ($ $
($1/60 , 1$ )$ , 1$ )$ , $ ($ $ ($1/30 , 1$ )$ , 1$ )$ , $ ($ $
($1/24 , 1$ )$ , 1$ )$ , $ ($ $ ($1/20 , 1$ )$ , 1$ )$ , $ ($ $
($1/12 , 1$ )$ , 4$ )$ , $ ($ $ ($1/10 ,  1$ )$ , 1$ )$ , $ ($ $
($7/60 , 1$ )$ , 1$ )$ , $ ($ $ ($1/8 , 1$ )$ , 1$ )$ , $ ($ $
($3/20  ,  1 $ )$ , 1$ )$ , $ ($ $ ($1/6 , 1$ )$ , 4$ )$ , $ ($ $
($11/60 , 1$ )$ , 1$ )$ , $ ($ $ ($13/60 , 1$ )$  , 1$ )$ , $ ($ $
($7/30 , 1$ )$ , 1$ )$ , $ ($ $ ($1/4 , 1$ )$ , 4$ )$ , $ ($ $
($7/24 , 1$ )$ , 1$ )$ , $ ($ $ ($3/10 , 1$ )$ , 1$ )$ , $ ($ $ ($
19/60  ,  1$ )$  ,  1$ )$ , $ ($ $ ($1/3 , 1$ )$ , 1$ )$  , $ ($ $
($1/3 , 0$ )$ , 1$ )$ , $ ($ $ ($ 11/30 , 1$ )$ , 1$ )$ , $ ($ $
($3/8 , 1$ )$ , 1$ )$ , $ ($ $ ($23/60 , 1$ )$ , 1$ )$ , $ ($ $
($5/12  , 1$ )$ , 1$ )$ , $ ($ $ ($13/30 , 1$ )$ , 1$ )$ ,  $ ($ $
($9/20 , 1$ )$ , 1$ )$ , $ ($ $ ($11/24 , 1$ )$ , 1$ )$ , $ ($ $
($1/2 , 1$ )$ , 1$ )$ , $ ($ $ ($1/2 , 0$ )$ , 1$ )$ , $ ($ $
($31/60 , 1$ )$ , 1$ )$ , $ ($ $ ($7/12 , 1$ )$ , 3$ )$ , $ ($ $
($ 13/20 , 1$ )$ , 1$ )$ , $ ($ $ ($2/3 , 0$ )$ , 1$ )$ , $ ($ $
($43/60 , 1$ )$ , 1$ )$ , $ ($ $ ($19/24 , 1$ )$ , 1$ ).$

Let us use the notation $\vec v_1=(0 ,12),$ $\vec v_2=(3 , 6),$
$\vec v_3=(6,4),$ $\vec v_4=(5,0),$ $\tau_1 = convex \;hull\{\vec
v_1, \vec v_2\},$ $\tau_2 = convex \;hull\{\vec v_2, \vec v_3\},$
$\tau_3 = convex \;hull\{\vec v_3, \vec v_4\}.$ Then we have
$$ supp(V_{\tau_1})= \bZ^2\cap \{ convex \;hull\{\vec v_1, \vec v_2,
\vec v_1 + \vec v_2\}^{int} \cup convex \;hull\{0, \vec
v_2\}^{int} \}.$$
$$ supp(V_{\tau_2})= \bZ^2\cap \{ convex \;hull\{\vec v_2, \vec v_3,
\vec v_2 + \vec v_3\}^{int}\cup convex \;hull\{0, \vec v_2\}^{int}
\cup convex \;hull\{0, \vec v_3\}^{int} \cup 2\{ 0\} \}.$$
$$ supp(V_{\tau_3})= \bZ^2\cap \{ convex \;hull\{\vec v_3, \vec v_4,
\vec v_3 + \vec v_4\}^{int}\cup convex \;hull\{0, \vec v_3\}^{int}
\}.$$ As we see there are repetitive appearances of $convex
\;hull\{0, \vec v_2\}^{int},$ $convex \;hull\{0, \vec v_3\}^{int}$
and $\{ 0\}$ each of them twice. Thus the summation $(1.5)$ must
be taken in the following way,
$$ A(f_1)\cong \bZ^2\cap \{ convex \;hull\{\vec v_1, \vec v_2,
\vec v_1 + \vec v_2\}^{int} \cup convex \;hull\{\vec v_2, \vec
v_3, \vec v_2 + \vec v_2\}^{int} \leqno(2.1.1)$$ $$\cup convex
\;hull\{\vec v_3, \vec v_4, \vec v_3 + \vec v_4\}^{int} \cup
convex \;hull\{0, 2\vec v_2\}^{int}\cup convex \;hull\{0, 2\vec
v_3\}^{int}\}. \leqno(2.1.2)$$ Here it is worthy to notice that
$$ convex \;hull\{0, 2\vec
v_2\}^{int} \cong convex \;hull\{0, \vec v_2\}^{int}\cup \{\vec
v_2\} \cup c_{\tau_1}(convex \;hull\{0, \vec
v_2\}^{int}),\leqno(2.1.3)$$
$$\cong convex \;hull\{0, \vec v_2\}^{int} \cup \{\vec v_2\} \cup
c_{\tau_2}(convex \;hull\{0, \vec v_2\}^{int}),$$
$$ convex \;hull\{0, 2\vec
v_3\}^{int} \cong convex \;hull\{0, \vec v_3\}^{int} \cup \{\vec
v_3\} \cup c_{\tau_2}(convex \;hull\{0, \vec v_3\}^{int}),$$
$$\cong convex \;hull\{0, \vec v_3\}^{int} \cup \{\vec v_3\} \cup
c_{\tau_3}(convex \;hull\{0, \vec v_3\}^{int}).$$ Douai
\cite{Douai1} exploits this kind of ambiguities scrupulously to
calculate the Gauss-Manin system. We can calculate by hands the
spectral pairs above in evaluating the monomials $[i], 1 \leq i
\leq 94$ modulo Jacobian ideal of $f_1$ by means of a piecewise
linear function,
$$
\begin {array}{ccccccc}
h(i,j)& =& \frac{i}{6} + \frac{j}{12} & for \;(i,j)\in \bar
B_{\tau_1}
\\
 &=& \frac{i}{12} + \frac{j}{8}
& for \;(i,j)\in \bar B_{\tau_2}\\
&=& \frac{i}{15} + \frac{3j}{20}
& for \;(i,j)\in \bar B_{\tau_3},\\
\end{array}$$
according to their classification into $\bar B_{\tau_1},$ $\bar
B_{\tau_2},$ $\bar B_{\tau_3}$ (closures of parallelepipeds
introduced in (1.17)).

 For example $$h([69]')-1 = \frac{3}{12} +
\frac{2}{8}-1 =\frac{3}{6} + \frac{2}{12}-1= -\frac{1} {2}, $$
which gives the spectral pair $((-\frac{1} {2},2),1).$ Here the
weight filtration index $2$ indicates that $[69]' \in cone (\tau_2
\cap \tau_3).$ In a similar way
$$h([71]')-1 = \frac{1}{6} + \frac{2}{12}-1
=\frac{1}{12} + \frac{2}{8}-1= -\frac{2} {3}, $$ that gives the
spectral pair  $((-\frac{2} {3},2),1).$

 {\bf 2.2} Let us treat the case studied by
\cite{Dan} as an example.  $$f_2=x4+y4+z8+x2z2+y2z2.$$
 The
ring $A(f_2)$ with rank $31$ has the following basis,

\noindent
 [1]$=x4y2 $ , [2]$=xy2z3$ , [3]$=z8 $ , [4]$=z7 $ ,
[5]$ $=z6 $ , [6]$=z5 $ , [7]$=yz4 $ , [8]$=xz4 $ , [9]$=z4 $ ,
[10]$=y3z2 $ , [11]$=x3y2 $ , [12]$=z3 $ , [13]$=y2z2 $ ,
[14]$=xyz2 $ , [15]$=yz2 $ , [16]$=xz2 $ , [17]$=z2 $ , [18]$=y2z
$ , [19]$=xyz $ , [20]$=yz $ , [21]$=xz $ , [22]$=z $ , [23]$=x2y2
$ , [24]$=xy2 $ , [25]$=y2 $ , [26]$=x2y $ , [27]$=xy $ , [28]$=y
$ , [29]$=x2 $ , [30]$=x $ , [31]$ $=1.$

This result is slightly different from what  $\sf SINGULAR$ gives
us due to the reason mentioned in Remark ~\ref{remark61}.

$\sf SINGULAR$ calculates the spectral pairs as follows,

\noindent $((-1/4 , 2$ )$ , 1$ )$ , ((0 , 3$ )$ , 1$ )$ , ((0 , 2$
)$ , 2$ )$ , ((1/8 , 2$ )$ , 1$ )$ , ((1/4 , 2$ )$ , 6$ )$ , ((3/8
, 2$ )$ , 1$ )$ , ((1/2 , 2$ )$ , 7$ )$  , ((5/8 , 2$ )$ , 1$ )$ ,
((3/4 , 2$ )$ , 6$ )$ , ((7/8 , 2$ )$ , 1$ )$ , ((1 , 2$ )$ , 2$
)$ , ((1 , 1$ )$ , 1$ )$ , ((5/4 , 2$ )$ , 1$ )$.$

Let us denote by $\Gamma_1$ the convex hull of $\{(0,0,0),(0,0,8),
(2,0,2), (0,2, 2) \},$ $\Gamma_2$ that of $\{(0, 0, 0),$ $(2, 0,
2),$ $(0, 2, 2), (0, 4, 0) \},$ $\Gamma_3$ that of
$\{(0,0,0),(2,0,2), (4,0,0), (0,4, 0) \}.$ The piecewise linear
function $h(i_1,i_2,i_3)$ is given by the following,
$$
\begin {array}{ccccccc}
h(i_1,i_2,i_3)& =& \frac{ i_1+i_2+i_3}{4} & for \;
(i_1,i_2,i_3)\in \bar B_{\Gamma_2}\cup \bar B_{\Gamma_3}
\\
 &=& \frac{ 3(i_1+i_2)+i_3}{8}
& for \;(i_1,i_2,i_3)\in \bar B_{\Gamma_1}. \\
\end{array}$$
We remark that $[1]'=x5y3z \in F^0/F^1(A_{\Gamma_1}) , $ while
$[31]'=xyz \in F^2(A_{\Gamma_1}).$ The point $[1]'$ is the
canonical copy of $[31]'$ with respect to $\Gamma_1.$ The point
$[2]' \in H^{1,1}_{\chi=1}$ is located on $cone (\Gamma_1 \cap
\Gamma_2)$ with spectral pair $(1 , 1)$ which is the canonical
copy of $[22]'\in H^{2,2}_{\chi=1}$ with respect to $\Gamma_1$
whose spectral pair is $(0 , 3).$ The points $[10]' , $ $[11]' \in
H^{1,2}_{\chi=1}$ are located on the 2-dimensional open skeleton
of $2\hat \Gamma_1$ and they give spectral pairs $((1 , 2) , 2).$
They are the canonical copies of $[28]',$ $[30]' \in
H^{2,1}_{\chi=1}$ with spectral pairs $((0,2),2).$ All other
integer points are located in the interior of $cone(\Gamma_1)\cup
cone(\Gamma_1 \cup \Gamma_2)$ with weight filtration $w=2$ and
they correspond to $H^{1,1}_{\chi\not=1} \oplus
H^{2,0}_{\chi\not=1}\oplus H^{0,2}_{\chi\not=1}.$ Here we recall
that $\chi = e^{-2 \pi h([i]') \sqrt -1}$ for each basis element
$[i].$

{\bf 2.3} Next we consider the case that B.Malgrange (in a letter
to the editor of Inventiones Mathematicae) used to demonstrate
that a maximum size Jordan cell (=the dimension $n$) of the
monodromy $T$ (or equivalently that of $T_u$) really appears,
$$f_3=x8+y8+z8+x2y2z2.$$
The MHS and essentially the spectral pairs of $Spp(f_3)$ are
described in detail in \cite {Steen1}, (3.15).

The ring $A(f_3)$ with rank 215 has the following basis,

 \noindent $\;$[$1$]$=z16 $, [$2$]$=z15 $, [$3$]$=z14 $,
[$4$]$=z13 $, [$5$]$=z12 $, [$6$]$=z11 $, [$7$]$=z10 $,
[$8$]$=xyz9 $, [$9$]$=yz9 $, [$10$]$=xz9 $, [$11$]$=z9 $,
[$12$]$=y7z8 $, [$13$]$=y6z8 $, [$14$]$=y5z8 $, [$15$]$=x12$,
[$16$]$=y3z8 $, [$17$]$=y2z8 $, [$18$]$=xyz8 $, [$19$]$=yz8 $,
[$20$]$=x7z8 $, [$21$]$=x6z8 $, [$22$]$=x5z8 $, [$23$]$=y12 $,
[$24$]$=x3z8 $, [$25$]$=x2z8 $, [$26$]$=xz8 $, [$27$]$=z8 $,
[$28$]$=y7z7 $, [$29$]$=y6z7 $, [$30$]$=y5z7 $, [$31$]$=y4z7 $,
[$32$]$=y3z7 $, [$33$]$=y2z7 $, [$34$]$=xyz7 $, [$35$]$=yz7 $,
[$36$]$=x7z7 $, [$37$]$=x6z7 $, [$38$]$=x5z7 $, [$39$]$=x4z7 $,
[$40$]$=x3z7 $, [$41$]$=x2z7 $, [$42$]$=xz7 $, [$43$]$=z7 $,
[$44$]$=y7z6 $, [$45$]$=y6z6 $, [$46$]$=y5z6 $, [$47$]$=y4z6 $,
[$48$]$=y3z6 $, [$49$]$=y2z6 $, [$50$]$=xyz6 $, [$51$]$=yz6 $,
[$52$]$=x7z6 $, [$53$]$=x6z6 $, [$54$]$=x5z6 $, [$55$]$=x4z6 $,
[$56$]$=x3z6 $, [$57$]$=x2z6 $, [$58$]$=xz6 $, [$59$]$=z6 $,
[$60$]$=y7z5 $, [$61$]$=y6z5 $, [$62$]$=y5z5 $, [$63$]$=y4z5 $,
[$64$]$=y3z5 $, [$65$]$=y2z5 $, [$66$]$=xyz5 $, [$67$]$=yz5 $,
[$68$]$=x7z5 $, [$69$]$=x6z5 $, [$70$]$=x5z5 $, [$71$]$=x4z5 $,
[$72$]$=x3z5 $, [$73$]$=x2z5 $, [$74$]$=xz5 $, [$75$]$=z5 $,
[$76$]$=y7z4 $, [$77$]$=y6z4 $, [$78$]$=y5z4 $, [$79$]$=y4z4 $,
[$80$]$=y3z4 $, [$81$]$=y2z4 $, [$82$]$=xyz4 $, [$83$]$=yz4 $,
[$84$]$=x7z4 $, [$85$]$=x6z4 $, [$86$]$=x5z4 $, [$87$]$=x4z4 $,
[$88$]$=x3z4 $, [$89$]$=x2z4 $, [$90$]$=xz4 $, [$91$]$=z4 $,
[$92$]$=y7z3 $, [$93$]$=y6z3 $, [$94$]$=y5z3 $, [$95$]$=y4z3 $,
[$96$]$=y3z3 $, [$97$]$=y2z3 $, [$98$]$=xyz3 $, [$99$]$=yz3 $,
[$100$]$=x7z3 $, [$101$]$=x6z3 $, [$102$]$=x5z3 $, [$103$]$=x4z3
$, [$104$]$=x3z3 $, [$105$]$=x2z3 $, [$106$]$=xz3 $, [$107$]$=z3
$, [$108$]$=y7z2 $, [$109$]$=y6z2 $, [$110$]$=y5z2 $,
[$111$]$=y4z2 $, [$112$]$=y3z2 $, [$113$]$=y2z2 $, [$114$]$=xyz2
$, [$115$]$=yz2 $, [$116$]$=x7z2 $, [$117$]$=x6z2 $, [$118$]$=x5z2
$, [$119$]$=x4z2 $, [$120$]$=x3z2 $, [$121$]$=x2z2 $, [$122$]$=xz2
$, [$123$]$=z2 $, [$124$]$=xy7z $, [$125$]$=y7z $, [$126$]$=xy6z
$, [$127$]$=y6z $, [$128$]$=xy5z $, [$129$]$=y5z $, [$130$]$=xy4z
$, [$131$]$=y4z $, [$132$]$=xy3z $, [$133$]$=y3z $, [$134$]$=xy2z
$, [$135$]$=y2z $, [$136$]$=x7yz $, [$137$]$=x6yz $, [$138$]$=x5yz
$, [$139$]$=x4yz $, [$140$]$=x3yz $, [$141$]$=x2yz $, [$142$]$=xyz
$, [$143$]$=yz $, [$144$]$=x7z $, [$145$]$=x6z $, [$146$]$=x5z $,
[$147$]$=x4z $, [$148$]$=x3z $, [$149$]$=x2z $, [$150$]$=xz $,
[$151$]$=z $, [$152$]$=x7y7 $, [$153$]$=x6y7 $, [$154$]$=x5y7 $,
[$155$]$=x4y7 $, [$156$]$=x3y7 $, [$157$]$=x2y7 $, [$158$]$=xy7 $,
[$159$]$=y7 $, [$160$]$=x7y6 $, [$161$]$=x6y6 $, [$162$]$=x5y6 $,
[$163$]$=x4y6 $, [$164$]$=x3y6 $, [$165$]$=x2y6 $, [$166$]$=xy6 $,
[$167$]$=y6 $, [$168$]$=x7y5 $, [$169$]$=x6y5 $, [$170$]$=x5y5 $,
[$171$]$=x4y5 $, [$172$]$=x3y5 $, [$173$]$=x2y5 $, [$174$]$=xy5 $,
[$175$]$=y5 $, [$176$]$=x7y4 $, [$177$]$=x6y4 $, [$178$]$=x5y4 $,
[$179$]$=x4y4 $, [$180$]$=x3y4 $, [$181$]$=x2y4 $, [$182$]$=xy4 $,
[$183$]$=y4 $, [$184$]$=x7y3 $, [$185$]$=x6y3 $, [$186$]$=x5y3 $,
[$187$]$=x4y3 $, [$188$]$=x3y3 $, [$189$]$=x2y3 $, [$190$]$=xy3 $,
[$191$]$=y3 $, [$192$]$=x7y2 $, [$193$]$=x6y2 $, [$194$]$=x5y2 $,
[$195$]$=x4y2 $, [$196$]$=x3y2 $, [$197$]$=x2y2 $, [$198$]$=xy2 $,
[$199$]$=y2 $, [$200$]$=x7y $, [$201$]$=x6y $, [$202$]$=x5y $,
[$203$]$=x4y $, [$204$]$=x3y $, [$205$]$=x2y $, [$206$]$=xy $,
[$207$]$=y $, [$208$]$=x7 $, [$209$]$=x6 $, [$210$]$=x5 $,
[$211$]$=x4 $, [$212$]$=x3 $, [$213$]$=x2 $, [$214$]$=x $,
[$215$]$=1.$ Here we have chosen $[15],$ $ [23]$ differently from
the computation given by $\sf SINGULAR$ due to the reason
mentioned in Remark ~\ref{remark61}.

The spectral pairs are calculated by $\sf SINGULAR$ as follows,
\noindent $($$($-1/2 , 4) , 1) , $($ $($-3/8 , 3) , 3) , $($
$($-1/4 , 3) , 3) , $($ $($-1/4 , 2) , 3) , $($ $($-1/8 , 3) , 3)
, $($ $($-1/8 , 2) , 6) , $($ $($ 0 , 3) , 4) , $($ $($0 , 2) , 9)
, $($ $($1/8 , 3) , 3) , $($ $($1/8 , 2) , 15) , $($ $($1/4 , 3) ,
3) , $($ $($1/4 , 2) , 18) , $($ $($3/8 , 3) , 3) ,  $($ $($3/8 ,
2) , 21) , $($ $($1/2 , 2) , 25) , $($ $($5/8 , 2) , 21) , $($
$($5/8 , 1) , 3) , $($ $($3/4 , 2) , 18) , $($ $($3/4 , 1) , 3) ,
$($ $($7/ 8 , 2) , 15) , $($ $($7/8 , 1) , 3) , $($ $($1 , 2) , 9)
, $($ $($1 , 1) , 4) , $($ $($9/8 , 2) , 6) , $($ $($9/8 , 1) , 3)
, $($ $($5/4 , 2) , 3) , $($ $($5 /4 , 1) , 3) , $($ $($11/8 , 1)
, 3) , $($ $($3/2 , 0$ )$ , 1$ ).$

Let us denote by $\Gamma_1$ the convex hull of $\{\vec
v_3=(0,0,8), \vec v_0=(2,2,2), \vec v_1= (8,0,0) \},$ $\Gamma_2$
that of $\{ \vec v_0,$ $\vec v_1,$ $\vec v_2=(0, 8, 0) \},$
$\Gamma_3$ that of $\{\vec v_0, \vec v_2, \vec v_3 \}.$ The most
interesting monomials with spectral pairs $(0 , 3) ,  (0 , 2)$ are
the following
$$ [91] ,  [142] ,  [183] ,  [211] \in H^{2 , 2}_{\chi=1} \;\; with\; spp(f_3)=(0 , 3) , $$
$$ [99] ,  [106] ,  [112] ,  [121], [132], [148], [190], [197], [204]
\in  H^{2,1}_{\chi=1}\;\; with\; spp(f_3)=(0,2).$$ Monomials of
$H^{1,2}_{\chi=1}$ (with spectral pairs ((0,2),9)) are obtained as
the canonical copies of $H^{2,1}_{\chi=1}$ ( with spectral pairs
((1,2),9) ) with respect to properly chosen  2-faces. Namely,
$$ [30] ,  [38] ,  [45] ,  [53], [60], [68], [154], [161], [168]
\in  H^{1,2}_{\chi=1}\;\; with\; spp(f_3)=(1,2).$$ We see also,
$$ [5] ,  [8] ,  [15] ,  [23]
\in  H^{1,1}_{\chi=1}\;\; with\; spp(f_3)=(1,1).$$ We see also $
H^{2,2}_{\chi=-1} =\{[215]\}$ and $ H^{0,0}_{\chi=-1} =\{[1]\}.$
All other monomials are located in $B_{\Gamma_1}^{int} \cup$
$B_{\Gamma_2}^{int} \cup$ $B_{\Gamma_3}^{int}.$


\vspace{\fill}

%

\noindent

\begin{flushleft}
\begin{minipage} [t]{6.2cm}
  \begin{center}
{\footnotesize Independent University of Moscow\\
Bol'shoj Vlasievskij pereulok 11,\\
 Moscow, 121002,\\
Russia\\
{\it E-mails}:  tanabe@mccme.ru, tanabe@mpim-bonn.mpg.de \\}
\end{center}
\end{minipage}
\end{flushleft}

\end{document}